\documentclass[titlepage,twoside,12pt]{article}
\usepackage{amsmath}
\usepackage{amssymb}
\usepackage{amsfonts}
\begin{document}

{\LARGE \bf Surprising Properties of Non-Archimedean \\ \\ Field Extensions of the \\ \\ Real Numbers} \\ \\

{\bf Elem\'{e}r E ~Rosinger} \\ \\
{\small \it Department of Mathematics \\ and Applied Mathematics} \\
{\small \it University of Pretoria} \\
{\small \it Pretoria} \\
{\small \it 0002 South Africa} \\
{\small \it eerosinger@hotmail.com} \\ \\

\hfill {\it Dedicated to Marie-Louise Nykamp} \\ \\

{\bf Abstract} \\

This, under the present form, is a replacement that is a two part paper in which the new second part was brought together with my recently posted arxiv paper, upon the suggestion of the arxiv moderators. \\ \\

{\LARGE \bf PART I : \\ \\ 
Real Scalars with a Rich Structure \\ \\ of ''Walkable Worlds'' : \\ \\ Universes next to Universes, \\ \\ and/or within Universes ... \\ \\ and so on, ad infinitum ...} \\ \\

{\bf Abstract} \\

It is a rather universal tacit and unquestioned belief - and even more so among physicists - that there is one and only one set of real scalars, namely, the one given by the usual field $\mathbb{R}$ of real numbers, with its usual linear order structure on the geometric line. Such a dramatically limiting and thus harmful belief comes,
unknown to equally many, from the similarly tacit acceptance of the ancient Archimedean Axiom in Euclid's Geometry. The consequence of that belief is a similar belief in the uniqueness of the field $\mathbb{C}$ of complex numbers, and therefore, of the various spaces, manifolds, etc., be they finite or infinite dimensional, constructed upon the real or complex numbers, including the Hilbert spaces used in Quantum Mechanics.
An near total lack of awareness follows about the {\it rich self-similar} structure of various linearly ordered scalar fields obtained through the
{\it ultrapower} construction which extend the usual field $\mathbb{R}$ of real numbers. Such {\it ultrapower field} extensions contain as a rather small subset the usual field $\mathbb{R}$ of real numbers. The rich self-similar
structure of such field extensions is due to {\it infinitesimals}, and thus also of {\it infinitely large} elements in such fields, which make such fields {\it non-Archimedean}. With the concept of {\it walkable world}, which has highly intuitive
and pragmatic algebraic and geometric meaning, the mentioned rich self-similar structure is illustrated. The ultrapower fields presented can have a wide ranging relevance in Physics, among others, for a proper treatment of what
are usually called the "infinities in Physics". The ultrapower construction which gives such non-Archimedean fields is rather simple and elementary,
requiring only 101 Algebra. \\ \\

{\bf 0. Preliminaries} \\

It is a rather universal tacit and hardly ever questioned belief - and even more so among physicists - that there is {\it no}, and thus there {\it cannot} be any other geometric
infinite straight line, except the one given by the usual field $\mathbb{R}$ of real numbers, with its usual linear order structure. Such a belief comes,
unknown to equally many, from the similarly tacit acceptance of the ancient Archimedean Axiom in Euclid's Geometry. \\

On the other hand, as is well known in Mathematics, there are a variety of {\it far larger} infinite straight lines which have a {\it far more rich} structure, among
the better known ones being those given in Topology by the so called {\it long line}, or in Analysis by the field $^*\mathbb{R}$ of {\it nonstandard} real numbers. \\

Needless to say, restricting oneself to scalars given by the usual field $\mathbb{R}$ of real numbers, as well as those built upon them, like the usual field $\mathbb{C}$ of complex numbers, can lead to dramatic effects in thinking in Physics. In particular, it can lead to inappropriate mathematical models of already known concepts, phenomena, processes, etc., or worse yet, it can simply prevent the very emergence of concepts in Physics due to the sheer impossibility of their mathematical modelling based on the dramatically limiting usual field $\mathbb{R}$ of real numbers used since ancient times and till the present day. \\
The so called "infinities in Physics" are an example of inappropriate mathematical models, while regarding the concepts in Physics not yet formulated they are but anybody's guess ... \\
One possible such concept still missing in Physics, and with likely fundamental importance, is that of infinitely many different {\it levels of precision}, [13, subsection 3.3]. \\

In this paper, we shall focus on an large class of extensions of the usual field $\mathbb{R}$ of real numbers and of its usual linear order structure with its the geometric line, a class other than the mentioned long line or nonstandard reals. These extensions will be constructed in a surprisingly easy and elementary way, much unlike for instance, the long line or nonstandard real numbers. The respective construction is called the {\it ultrapower} construction, and the resulting linearly ordered field extensions of the usual field $\mathbb{R}$ of real numbers are called {\it ultrapower fields}. \\
The natural aspect of the ultrapower construction is in the fact that, as seen in Model Theory, which is a branch of Mathematical Logic, this construction is in fact often present in a large variety of important mathematical contexts. And the general lack of familiarity with that fact is only due to the reduced familiarity with Mathematical Logic, let alone, Model Theory, among mathematicians, and even more so, among physicists. \\

The importance of the existence and easy, elementary construction of such ultrapower fields is that, contrary to the general view in Physics, the real or complex scalars, as well as the variety of finite or infinite dimensional spaces or manifolds constructed upon them, including Hilbert spaces, can in fact have a surprisingly {\it rich self-similar} structure, one that is so far just about completely unknown in Physics, and one that can offer the modelling of a wide range of yet unknown phenomena, as well as the more appropriate modelling of the already known ones. \\

As mentioned, the consequences of the present highly limiting belief regarding the structure of the geometric straight line, and thus of all the scalars and spaces built upon it, can lead to unfortunate effects in Physics. Among them may simply be the inability to model properly important physical phenomena, and here one can immediately think
about the considerable troubles the so called "infinities in Physics" do cause, among others, in Quantum Field Theory. Yet more grave consequences could
also come from the sheer impossibility to conceive of, let alone model properly important physical concepts, phenomena, processes, etc. \\

Such and relates issue have been addressed in certain of their more general features in [2-14]. \\

In [11-13] it has been shown that certain most basic results in Relativity, Quantum Mechanics and Quantum Computation can in an easy and natural manner be reformulated in terms of the far more rich space-time structures built upon the general ultrapower construction of which the ultrapower fields are but a particular case. \\

Here, with the help of a highly intuitive and pragmatic algebraic and geometric concept, namely,
that of a {\it walkable world}, the surprisingly {\it rich self-similar} structure of a large class of infinite geometric lines is presented, lines which correspond to linearly ordered ultrapower fields containing the usual field $\mathbb{R}$ of real numbers. \\

Since the Arithmetisation of Geometry by Descartes, the usual infinite geometric straight line is associated in the well known manner with the field
$\mathbb{R}$ of the usual real numbers. \\

Such an association, fortunately, can easily be extended to hold between various far larger infinite geometric straight lines, lines with a rich self-similar
structure, and on the other hand, linearly ordered fields which contain the field $\mathbb{R}$ of the usual real numbers as a rather small subset. \\
In view of that association, an association which will easily follow from the constructions in the sequel, we shall only deal with the linearly ordered ultrapower
field extensions of the field $\mathbb{R}$ of the usual real numbers. \\

These ultrapower field extensions will be obtained as particular cases of the so called {\it ultrapower} construction. \\ \\

{\bf 1. Linearly Ordered Fields as Ultrapowers} \\

There is a simple and elementary way to construct linearly ordered fields which contain as a rather small subset the usual field $\mathbb{R}$ of the real numbers. This
construction, called the {\it ultrapower} construction, proceeds as follows, [2-14]. \\

Given an infinite index set $\Lambda$, we take on it any {\it ultrafilter} ${\cal U}$, see Appendix, which satisfies the condition \\

(1.1)~~~ $ {\cal F}_{re} ( \Lambda ) \subseteq {\cal U} $ \\

where \\

(1.2)~~~ $ {\cal F}_{re} ( \Lambda ) = \{~ I \subseteq \Lambda ~~|~~ \Lambda \setminus I ~~\mbox{is finite} ~\} $ \\

is called the Frech\`{e}t filter on $\Lambda$. Further, we define on the set of all real valued functions defined on $\Lambda$, that is, on the set \\

$~~~~~~ \mathbb{R}^\Lambda = \{~ x : \Lambda \longrightarrow \mathbb{R} ~\} $ \\

the equivalence relation $\approx_{\cal U}$ by \\

(1.3)~~~ $ x \approx_{\cal U} y ~~\Longleftrightarrow~~
                        \{~ \lambda \in \Lambda ~|~ x ( \lambda ) = y ( \lambda ) ~\} \in {\cal U} $ \\

Finally, through the usual set-theoretic quotient construction, we obtain the {\it ultrapower field} \\

(1.4)~~~ $ \mathbb{F}_{\cal U} = \mathbb{R}^\Lambda / \approx_{\cal U} $ \\

which proves to have the following two properties. \\

First, the mapping \\

(1.5)~~~ $ \mathbb{R} \ni r \longmapsto ( u_r )_{\cal U} \in \mathbb{F}_{\cal U} $ \\

is an {\it embedding of fields} in which $\mathbb{R}$ is a strict subset of $\mathbb{F}_{\cal U}$, where $u_r \in
\mathbb{R}^\Lambda$ is defined by $u_r ( \lambda ) = r$, for $\lambda \in \Lambda$, while $( u_r )_{\cal U}$ is the coset of $u_r$ with respect to the equivalence relation $\approx_{\cal U}$. For simplicity we shall denote
$( u_r )_{\cal U} = r$, for $r \in \mathbb{R}$, and thus (1.5) takes the form \\

(1.6)~~~ $ \mathbb{R} \ni r \longmapsto r \in \mathbb{F}_{\cal U},~~~
                         \mbox{or simply}~~~ \mathbb{R} \subsetneqq \mathbb{F}_{\cal U} $ \\

The field operations of addition and multiplication on $\mathbb{F}_{\cal U}$ are simply given by the following {\it
term-wise} operations. If $\xi = ( x )_{\cal U},~ \eta = ( y )_{\cal U} \in \mathbb{F}_{\cal U}$, where $x, y \in
\mathbb{R}^\Lambda$, then \\

(1.7)~~~ $ \xi + \eta = ( x + y )_{\cal U},~~~ \xi \eta = ( x y )_{\cal U} $ \\

Second, on $\mathbb{F}_{\cal U}$ we have the {\it linear order} defined by \\

(1.8)~~~ $ ( x )_{\cal U} \leq ( y )_{\cal U} ~~\Longleftrightarrow~~
             \{~ \lambda \in \Lambda ~|~ x ( \lambda ) \leq y ( \lambda ) ~\} \in {\cal U} $ \\

where $x, y \in \mathbb{R}^\Lambda$. And this linear order is compatible with the field structure of $\mathbb{F}_{\cal
U}$, which means that, for $\xi, \eta, \theta \in \mathbb{F}_{\cal U}$, we have \\

(1.9)~~~ $ \xi \leq \eta ~~~\Longrightarrow~~~ \xi + \theta \leq \eta + \theta $ \\

while for for $\xi, \eta, \theta \in \mathbb{F}_{\cal U},~ \theta \geq 0$, we have \\

(1.10)~~~ $ \xi \leq \eta ~~~\Longrightarrow~~~ \xi \theta \leq \eta \theta$ \\

We recall that the field $\mathbb{R}$ of usual real numbers is linearly ordered, giving what in Geometry is called the real line. And with the above linear order in (1.8), $\mathbb{R}$ becomes a linearly ordered subfield of the much larger linearly ordered field $\mathbb{F}_{\cal U}$. \\
In other words, from geometric point of view, $\mathbb{F}_{\cal U}$ is a much more rich self-similar real line, than the usual real line given by $\mathbb{R}$. \\

Lastly, it should be noted that the nonstandard reals $^*\mathbb{R}$ are a particular case of the above ultrapower fields (1.4), see Remark 1.1. at the end of this section. \\

Here, it is worth recalling the general case of the above construction, even if for the sake of simplicity it will not be considered in the present paper. Namely, the above ultrapower construction is a particular case of what is well known in Model Theory - a branch of Mathematical Logic - under the name of {\it reduced powers}. Details on that general construction can be found in [2-14]. \\
This general construction of {\it reduced power algebras} goes as follows, and again, it only requires 101 Algebra. \\

Let ${\cal F}$ be any filter on $\Lambda$ which satisfies \\

(1.11)~~~ $ {\cal F}_{re} ( \Lambda ) \subseteq {\cal F} $ \\

We define on $\mathbb{R}^\Lambda$ the corresponding equivalence relation $\approx_{\cal F}$ by \\

(1.12)~~~ $ x \approx_{\cal U} y ~~\Longleftrightarrow~~
                        \{~ \lambda \in \Lambda ~|~ x ( \lambda ) = y ( \lambda ) ~\} \in {\cal F} $ \\

Then, through the usual quotient construction, we obtain the {\it reduced power algebra} \\

(1.13)~~~ $ \mathbb{A}_{\cal F} = \mathbb{R}^\Lambda / \approx_{\cal F} $ \\

which has the following two properties. The mapping \\

(1.14)~~~ $ \mathbb{R} \ni r \longmapsto ( u_r )_{\cal F} \in \mathbb{A}_{\cal F} $ \\

is an {\it embedding of algebras} in which $\mathbb{R}$ is a strict subset of $\mathbb{A}_{\cal F}$, where $u_r \in
\mathbb{R}^\Lambda$ is defined by $u_r ( \lambda ) = r$, for $\lambda \in \Lambda$, while $( u_r )_{\cal F}$ is the coset
of $u_r$ with respect to the equivalence relation $\approx_{\cal F}$. Further, on $\mathbb{A}_{\cal F}$ we have the {\it
partial order} \\

(1.15)~~~ $ ( x )_{\cal F} \leq ( y )_{\cal F} ~~\Longleftrightarrow~~
             \{~ \lambda \in \Lambda ~|~ x ( \lambda ) \leq y ( \lambda ) ~\} \in {\cal F} $ \\

where $x, y \in \mathbb{R}^\Lambda$. \\ \\

{\bf Remark 1.1.} \\

It is worth recalling that Abraham Robinson's field $^*\mathbb{R}$ of nonstandard reals can be obtained by (1.4) above. However, what complicates considerably the resulting Nonstandard Analysis is the intent to have the so called Transfer Principle operating in it. And according to that principle every property of the field $\mathbb{R}$ of the usual real numbers, property which cam be formulated in terms of First Order Predicate Logic, should also hold for the field $^*\mathbb{R}$ of nonstandard reals. \\

However, in Nonstandard Analysis it appears that one is obliged to pay a high price for achieving that intent. Indeed, the resulting necessary technical edifice makes Nonstandard Analysis quite hard even for a large majority of mathematicians, let alone physicists. On the other hand, the advantage obtained from the Transfer Principle turns out to be rather limited, since by far most of the properties of interest of $\mathbb{R}$, and let alone of $^*\mathbb{R}$, simply cannot be formulated in terms of First Order Predicate Logic. \\

On the contrary, here, as well as in [2-14], we do {\it not} aim to obtain any kind of transfer principle. Instead, we only make use of the {\it reduced power} construction. Consequently, the technical part can be limited to the simplicity and ease of 101 Algebra. \\ \\

{\bf 2. Walkable Worlds ...} \\

An essential property of the ultrapower fields (1.4) is that they are {\it no} longer Archimedean. In other words, unlike
the usual field given by the real line $\mathbb{R}$, they do {\it not} satisfy the Archimedean Axiom \\

(2.1)~~~ $ \exists~~ u > 0 ~:~ \forall~~ v > 0 ~:~ \exists~~ n \in \mathbb{N} ~:~ n u > v $ \\

As it happens, and it is still seldom realized, the fact that the ultrapower fields (1.4) are {\it non-Archimedean}, and
also are {\it larger} than $\mathbb{R}$, gives them an extremely rich both {\it local} and {\it global} structure. And
this results in a surprising {\it involved self-similar} structure. \\

Here we shall illustrate that fact with the help of the concept of {\it walkable world}, a concept which is highly
intuitive in its pragmatic algebraic and geometric meaning. \\

Let us briefly recall some of the basic features of the ultrapower fields (1.4), features which give their rich
self-similar structure. \\

One that follows immediately from the fact that they are non-Archimedean is that their elements $t \in \mathbb{F}_{\cal U}$ are of {\it three} kind,
namely, {\it infinitesimal, finite}, and {\it infinitely large}, as defined by the following respective conditions \\

(2.2)~~~ $ \forall~~ r \in \mathbb{R},~ r > 0 ~:~ t \in ( - r, r ) $ \\

(2.3)~~~ $ \exists~~ r \in \mathbb{R},~ r > 0 ~:~ t \in ( - r, r ) $ \\

(2.4)~~~ $ \forall~~ r \in \mathbb{R},~ r > 0 ~:~ t \notin ( - r,  r ) $ \\

where for $a, b \in \mathbb{F}_{\cal U}$, we denote as usual $( a, b ) = \{ s \in \mathbb{F}_{\cal U} ~|~ a < s < b \}$. Now,
following Leibniz, one denotes \\

(2.5)~~~ $ monad ( 0 ) = \{~ t \in \mathbb{F}_{\cal U} ~|~~ t ~~\mbox{is~infinitesimal}~ \} $ \\

and calls it the {\it monad} of $0 \in \mathbb{F}_{\cal U}$, while following Keisler, [1], one denotes \\

(2.6)~~~ $ Gal ( 0 ) = \{~ t \in \mathbb{F}_{\cal U} ~|~~ t ~~\mbox{is~finite}~ \} $ \\

and calls it the {\it Galaxy} of $0 \in \mathbb{F}_{\cal U}$. \\

It is easy to see that \\

(2.7)~~~ $ Gal ( 0 ) = \bigcup_{r \in \mathbb{R}} monad ( r ) $ \\

where for $t \in \mathbb{F}_{\cal U}$, we denote \\

(2.8)~~~ $ monad ( t ) = t + monad ( 0 ) $ \\

Finally \\

(2.9)~~~ $ \mathbb{F}_{\cal U} \setminus Gal ( 0 ) $ \\

is the set of infinitely large elements in the ultrapower field $\mathbb{F}_{\cal U}$. \\

In this way, all the elements of $\mathbb{F}_{\cal U}$, be they infinitesimal, finite, or infinitely large, have been expressed respectively in (2.5) by the monad of $0 \in \mathbb{F}_{\cal U}$, in (2.7) by the Galaxy of $0 \in \mathbb{F}_{\cal U}$, and in (2.9). And as one notes, all these sets can in fact be expressed in terms of the monad of $0 \in \mathbb{F}_{\cal U}$ alone. \\

And now, to the {\it walkable worlds} ... \\

Let $t, u \in \mathbb{F}_{\cal U}$, with $u > 0$. Then we denote \\

(2.10)~~~ $ \begin{array}{l}
                WW ( t, u ) =
                  \{~ s \in \mathbb{F}_{\cal U} ~|~ \exists~ n \in \mathbb{N} ~:~ s \in ( t - n u, t + n u ) ~\} = \\ \\
                  ~~~~~~~~~~~~~ = \bigcup_{\, n \in \mathbb{N}}~ ( t - n u, t + n u )
            \end{array} $ \\

which is the set of elements $s \in \mathbb{F}_{\cal U}$ that can be reached in either direction starting at $t$, by a {\it finite} number of
steps of length $u$. Thus $WW ( t, u )$ is the {\it walkable world} around $t$, by steps of length $u$. \\

Obviously \\

(2.11)~~~ $ Gal ( 0 ) = WW ( 0, 1 ) $ \\

and \\

(2.12)~~~ $ s \in WW ( t, u ) ~~\Longrightarrow~~ WW ( s, u ) = WW ( t, u ) $ \\

while \\

(2.13)~~~ $ WW ( t, v ) = WW ( t, u ) $ \\

for all $v \in \mathbb{F}_{\cal U},~ v > 0$, such that \\

(2.14)~~~ $ m u \leq v \leq ( m + 1 ) u $ \\

for some $m \in \mathbb{N}$. \\

Furthermore, for every $t, u \in \mathbb{F}_{\cal U}$, with $u > 0$, we have the {\it order isomorphism} \\

(2.15)~~~ $ WW ( t, u ) \ni s \longmapsto ( s - t ) / u \in WW ( 0, 1 ) $ \\

thus any two walkable worlds are {\it order isomorphic}. \\

The non-Archimedean nature of $\mathbb{F}_{\cal U}$ results, among others, in \\

(2.16)~~~ $ monad ( 0 ) \neq WW ( t, u ),~~~ \mathbb{F}_{\cal U} \neq WW ( t, u ) $ \\

and in fact \\

(2.17)~~~ $ WW ( t, u ) $ is an infinitely small part of $ \mathbb{F}_{\cal U} $ \\

for all $t, u \in \mathbb{F}_{\cal U}$, with $u > 0$. \\ \\

{\bf Remark 2.1.} \\

In view of (2.11), in any given ultrapower field $\mathbb{F}_{\cal U}$, the whole of $Gal ( 0 )$ is only {\it one single}
walkable world, while as seen in the next section, there are {\it infinitely many} walkable worlds, either {\it disjoint} from one
another, or {\it nested} within one another. \\

Here is already one of the essential {\it differences} with ultrapower fields, when compared with the usual field of the real
line $\mathbb{R}$. Indeed, in the latter case, what corresponds to $Gal ( 0 ) = WW ( 0, 1 )$ is the {\it whole} of $\mathbb{R}$,
namely, for every $t, u \in \mathbb{R},~ u > 0$, we have \\

(2.18)~~~ $ \{~ s \in \mathbb{R} ~|~ \exists~ n \in \mathbb{N} ~:~ s \in ( t - n u, t + n u ) ~\} = \mathbb{R} $ \\

In other words, the usual field of the real line $\mathbb{R}$ is but only {\it one single} walkable world, while in the case of
ultrapower fields, each walkable world is merely an infinitely small part of such a field, as seen in (2.17). \\

This difference clearly illustrates how much limitation is imposed upon the structure of $\mathbb{R}$ by the acceptance of the Archimedean Axiom. \\ \\

{\bf 3. Universes Next To and/or Within Universes, and so on \\
        \hspace*{0.45cm} Ad Infinitum ...} \\

Let us see now the way two arbitrary walkable worlds $WW ( t, u )$ and $WW ( s, v )$, with $t, u, s, v \in
\mathbb{F}_{\cal U},~ u, v > 0$, can relate to one another. As it turns out, we can distinguish the following {\it three}
situations regarding their relative position to one another : \\

{\bf Case I} ~:~ The two walkable worlds are {\it the same}, namely \\

(3.1)~~~ $ WW ( t, u ) = WW ( s, v ) $ \\

{\bf Case II} ~:~ The two walkable worlds are {\it disjoint}, thus in particular, they are different, namely \\

(3.2)~~~ $ WW ( t, u ) \cap WW ( s, v ) = \phi $ \\

{\bf Case III} ~:~ The two walkable worlds are {\it not} disjoint, yet they are nevertheless {\it different}, namely \\

(3.3)~~~ $ WW ( t, u ) \cap WW ( s, v ) \neq \phi,~~~ WW ( t, u ) \neq WW ( s, v ) $ \\

Now, a specific feature with walkable worlds is that in the last case, we must necessarily have {\it either} the {\it nesting} \\

(3.3.1)~~~ $ WW ( t, u ) $ is an infinitesimal part of $ WW ( s, v ) $ \\

{\it or} the {\it nesting} \\

(3.3.2)~~~ $ WW ( s, v ) $ is an infinitesimal part of $ WW ( t, u ) $ \\

In other words, if two different walkable worlds intersect, that is, in the case of (3.3), then one of them must strictly
contain the other, and in fact, it must contain it as an infinitesimal part. \\

Furthermore, concerning (3.2), there are {\it infinitely} many walkable worlds which are pair-wise disjoint. As for (3.3), the
respective nestings in (3.3.1) and (3.3.2) always have an {\it infinite} depth, that is, there are always infinitely many
infinitely smaller and smaller walkable worlds nested in one another. \\

Let us go more into detail on the possible relation between two walkable worlds in the case (3.3). We note that, given $u, v
\in \mathbb{F}_{\cal U},~ u, v > 0$, we can in view of the linear order on $\mathbb{F}_{\cal U}$ always assume that \\

(3.4)~~~ $ v \leq u$ \\

Thus the following two alternatives result \\

(3.5)~~~ either $u / v$ is finite, or $u / v$ is infinitely large \\

and correspondingly, either \\

(3.6)~~~ $ \exists~~ n \in \mathbb {N},~ n \geq 1 :~ u \leq n v $ \\

or \\

(3.7)~~~ $ \forall~~ n \in \mathbb {N},~ n \geq 1 :~ n v \leq u $ \\

Now in view of (2.13), (2.14), the alternative (3.5) is not compatible with (3.3), since it leads to (3.1). Thus we remain with (3.7). And then (3.3.2) follows. \\

What is important to note here is that (3.3.2), with the assumption (3.4), will always have $u$ and $v$ in one and only one of the following {\it six} situations : \\

1) both $u$ and $v$ are infinitesimal \\

2) $u$ is finite, $v$ is infinitesimal \\

3) both $u$ and $v$ are finite \\

4) $u$ is infinitely large, $v$ is infinitesimal \\

5) $u$ is infinitely large, $v$ is finite \\

6) both $u$ and $v$ are infinitely large \\

We can conclude with \\

{\bf Theorem 3.1.} \\

The structure of the field $\mathbb{R}$ of usual real numbers in terms of walkable worlds is given by the simple relation \\

(3.8)~~~ $ \mathbb{R} = WW ( 0, 1 ) $ \\

On the other hand, the structure of the reduced power field $\mathbb{F}_{\cal U}$ in terms of walkable worlds is given by the relation \\

(3.9)~~~~ $ \mathbb{F}_{\cal U} = \bigcup_{s, u \in \, \mathbb{F}_{\cal U},\, u > \, \, 0}\, WW ( s, u ) $ \\

and the walkable worlds in(3.9) form a {\it nesting} on $\mathbb{F}_{\cal U}$ in the above richly structured manner.

\hfill $\Box$ \\

Here, related with Theorem 3.1. above, we recall \\

{\bf Definition 3.1.} \\

A family of nonvoid subsets $A_i \subseteq X$, with $i \in I$, is called a {\it nesting} on $X$, if and only if \\

(3.10)~~~ $ X = \bigcup_{i \in I} \, A_i $ \\

and, for $i, j \in I$, we have \\

(3.11)~~~ $ A_i \cap A_j \neq \phi ~~\Longrightarrow~~ A_i \subseteq A_j ~\mbox{or}~ A_j \subseteq A_i $ \\ \\

{\bf 4. The Rich Self-Similar Structure of Ultrapower Fields} \\

 For the sake of easier comparison, let us recall the {\it self-similarity} properties of the usual field $\mathbb{R}$ of
 real numbers. In this regard, we have the {\it self-similarity} property given by the following {\it bijective, order reversing} mapping \\

(4.1)~~~ $ \mathbb{R} \setminus ( -1, 1 ) \ni r \longmapsto 1 / r \in [ -1, 1 ] \setminus \{ 0 \} $ \\

thus the unbounded set \\

$~~~~~~ \mathbb{R} \setminus ( -1, 1 ) = ( -\infty, -1 ] \cup [ 1, \infty ) $ \\

has through the mapping (4.1) the inverse linear order structure of the bounded set \\

$~~~~~~ [ -1, 1 ] \setminus \{ 0 \} = [ -1, 0 ) \cup ( 0, 1 ] $ \\

Now by translation and scaling, we obtain the family of self-similarities of the usual field $\mathbb{R}$ of
 real numbers, given by the bijective, order reversing mappings \\

(4.2)~~~ $ \mathbb{R} \setminus ( -a, a ) \ni r \longmapsto
                         ( 1 / r ) + r_0 \in [ r_0 - \frac{1}{a}, r_0 + \frac{1}{a} ] \setminus \{ r_0 \} $ \\

where $r_0, a \in \mathbb{R}, a > 0$. \\

Here we can note that none of the self-similarities (4.2) refers to the structure at any given point $r_0 \in \mathbb{R}$, but only to the structure of the sets \\

(4.3)~~~ $ [ r_0 - a, r_0 + a ] \setminus \{ r_0 \} = [ r_0 - a, r_0 ) \cup ( r_0, r_0 + a ],~~~~ a > 0 $ \\

around points $r_0 \in \mathbb{R}$, sets which are whole neighbourhoods of $r_0$ from which, however, the point $r_0$ itself has been taken out. This is obviously inevitable, since each point $r_0 \in \mathbb{R}$ is at a finite strictly positive - thus not infinitesimal - distance from any other point in $\mathbb{R}$. \\

In addition, we also have the self-similarities \\

(4.4)~~~ $ \mathbb{R} \stackrel{f}\longrightarrow ( a, b ) $ \\

where $-\infty \leq a < b \leq \infty$, while $f$ can be any bijective order preserving continuous mapping. \\

On the other hand, with the ultrapower fields  $ \mathbb{F}_{\cal U}$, their self-similar structures are far more {\it
rich}, due to the presence of their infinitesimals, and thus as well, of their infinitely large elements. Indeed, this
time, the self-similarities can also refer to the whole {\it monad} of each point, except for the point itself. \\

Let us start with a self-similarity of any ultrapower field $ \mathbb{F}_{\cal U}$ which does {\it not} exist in the
case of the usual real line $\mathbb{R}$. Namely, it is easy to see that we have the order reversing bijective mapping \\

(4.5)~~~ $ ( \mathbb{F}_{\cal U} \setminus Gal(0) ) \ni t \longmapsto 1 / t \in ( monad(0) \setminus \{ 0 \} ) $ \\

which means that the set of all infinitely large elements in $\mathbb{F}_{\cal U}$ has the inverse order structure of the set of infinitesimal elements from which one excludes 0. \\

This shows the {\it important} fact that the infinitesimally local structure, and on the other hand, the global structure of $\mathbb{F}_{\cal U}$ do in fact {\it mirror} one another, a property which has no correspondence in the case of the usual field $\mathbb{R}$ of real numbers. \\

Also, through translation and scaling, we have, for each $t_0, u \in \mathbb{F}_{\cal U}, u > 0$, the order reversing bijective mapping \\

(4.6)~~~ $  \mathbb{F}_{\cal U} \setminus ( -u, u ) \ni t \longmapsto
                      ( 1 / t ) + t_0 \in [ t_0 - \frac{1}{u}, t_0 + \frac{1}{u} ] \setminus \{ t_0 \} $ \\

where $\mathbb{F}_{\cal U} \setminus ( -u, u )$ will always contain infinitely large elements. \\

These again are self-similarities not present in the case of the usual real line $\mathbb{R}$. \\

Furthermore, in (4.6) we have a far more rich possibility for translations and scalings than in the usual case of the real line $\mathbb{R}$. Indeed, in addition to translations and scalings with non-zero finite elements $r_0, a \in \mathbb{R}, a > 0$, as in (4.3), we can now also translate and scale with all $t_0, u \in \mathbb{F}_{\cal U}, u > 0$, thus with all infinitely large elements, as well as with all infinitesimal elements, except for scaling with $0 \in \mathbb{F}_{\cal U}$. \\

Let us consider the above in some detail by listing the different possibilities for the sets \\

(4.7)~~~ $ [ t_0 - \frac{1}{u}, t_0 + \frac{1}{u} ] \setminus \{ t_0 \} $ \\

in (4.6). \\

First of all, these sets are no longer mere subsets in $\mathbb{R}$, but instead, they are subsets in $\mathbb{F}_{\cal U}$, and will always contain infinitesimals, since they contain nonvoid intervals. Furthermore, as seen below, they may also contain infinitely large elements. \\

Also, $t_0, u \in \mathbb{F}_{\cal U}, u > 0$ in (4.7) can independently be finite, infinitesimal, or infinitely
large, thus resulting in 9 possible combinations and 6 distinct outcomes regarding the set (4.7), which we list below. This is in sharp contradistinction with the case in (4.3) which applies to the real line $\mathbb{R}$. Indeed : \\

1) Let us start the listing of these 9 different cases and 6 distinct outcomes with both $t_0$ and $u$ being finite. Then obviously (4.7) is a subset of $Gal(0)$, and it has the finite, non-infinitesimal length $2u$. \\

2) When $t_0$ is finite and $u$ is infinitesimal, then the set (4.7) is infinitely large, and is no longer contained in $Gal ( 0 )$, however, it contains $Gal ( 0 ) \setminus \{ t_0 \}$. \\

3) If $t_0$ is finite, but $u$ is infinitely large, then (4.7) is again a subset of $Gal ( 0 )$, and in fact, it has the infinitesimal length $2 u$, which means that it is a subset of $monad ( t_0 )$. \\

4) Let us now assume that $t_0$ is infinitesimal and $u$ finite. Then regarding the set (4.1), we are back to case 1) above. \\

5) If both $t_0$ and $u$ are infinitesimal then the set (4.7) is as in 2) above. \\

6) When $t_0$ is infinitesimal and $u$ is infinitely large, the set (4.7) is as in 3) above. \\

7) Let us now take $t_0$ infinitely large and $u$ finite. Then the set (4.7) is disjoint from $Gal ( 0 )$, and it has the finite, non-infinitesimal length $2 u$. \\

8) When $t_0$ infinitely large and $u$ infinitesimal, then the set (4.7) is again not contained in $Gal ( 0 )$, and it has the infinitely large length $2 u$. Furthermore, depending on the relationship between $|t_0|$ and $1 / u$, it may, or it may not intersect $Gal ( 0 )$. \\

9) Finally, when both $t_0$ and $u$ are infinitely large, then the set (4.7) is disjoint from $Gal ( 0 )$, and it has the infinitesimal length $2 u$. \\

We conclude that he {\it local} structure of $\mathbb{F}_{\cal U}$ is given by \\

(4.8)~~~ $ Gal ( 0 ) = \bigcup_{r \in \mathbb{R}} \, ( r + monad ( 0 ) ) $ \\

while the {\it global} structure of $\mathbb{F}_{\cal U}$ is given by \\

(4.9)~~~ $ \mathbb{F}_{\cal U} = (~ \bigcup_{\lambda \in \Lambda} \, Gal ( -s_\lambda ) ~) ~\bigcup~ Gal ( 0 )
                                              ~\bigcup~ (~ \bigcup_{\lambda \in \Lambda} \, Gal ( s_\lambda ) ~) $ \\

where $\Lambda$ is an uncountable set of indices, while $s_\lambda \in \mathbb{F}_{\cal U}$ are positive infinite, and such that $s_\mu - s_\lambda$ is infinite, for $\lambda, \mu \in \Lambda, \lambda \neq \mu$. \\

Here we can point to a self-similar aspect of the interrelation between the local and global structure of $\mathbb{F}_{\cal U}$ which may remind us of a typical feature of {\it fractals}. Indeed, similar with (4.8), the relation (4.9) can also be expressed in terms monads, namely \\

(4.10)~~~ $ \mathbb{F}_{\cal U} = (~ \bigcup_{r \in \mathbb{R}, \,\, \lambda \in \Lambda} \,
                           (r - s_\lambda + monad ( 0 )) ~) ~\bigcup~ $ \\ \\
          \hspace*{3cm} $ (~ \bigcup_{r \in \mathbb{R}} \, ( r + monad ( 0 ) ) ~) $ \\ \\
          \hspace*{3cm} $~\bigcup~ (~ \bigcup_{r \in \mathbb{R}, \,\, \lambda \in \Lambda} \, ( r + s_\lambda +
                            monad ( 0 )) ~) $ \\

In this way, in view of (4.5), we obtain the self-similar order reversing bijection, which is now expressed solely in terms of $mon ( 0 )$, namely \\

(4.11)~~~ $ [~ (~ \bigcup_{r \in \mathbb{R}, \,\, \lambda \in \Lambda} \,
                       (r - s_\lambda + monad ( 0 )) ~) $ \\ \\
                        \hspace*{2cm} $   ~\bigcup~ (~ \bigcup_{r \in \mathbb{R}, \,\, \lambda \in \Lambda} \, ( r + s_\lambda + monad ( 0 )) ~) ~] \ni t \longmapsto $ \\ \\
                        \hspace*{7cm} $ \longmapsto 1 / t \in [~ monad  ( 0 ) \setminus \{ 0 \} ~] $ \\

and conversely \\

(4.12)~~~ $ [~ monad  ( 0 ) \setminus \{ 0 \} ~] \ni t \longmapsto $ \\ \\
               \hspace*{5cm} $ \longmapsto 1/ t \in [~ (~ \bigcup_{r \in \mathbb{R}, \,\, \lambda \in \Lambda} \,
                       (r - s_\lambda + monad ( 0 )) ~) $ \\ \\
                        \hspace*{7cm} $   ~\bigcup~ (~ \bigcup_{r \in \mathbb{R}, \,\, \lambda \in \Lambda} \, ( r + s_\lambda + monad ( 0 )) ~) ~] $ \\

As we can note, the above bijections in (4.11), (4.12) are given by the very simple algebraic, explicit, and order reversing mapping $s \longmapsto 1 / s$, which involves what is essentially a {\it field} operation, namely, division.
And these two bijections take the place of the much simpler order reversing bijections in the case of the usual real line $\mathbb{R}$, namely \\

(4.13)~~~ $ ( \mathbb{R} \setminus ( -1, 1) ) \ni r \longmapsto 1 / r \in  ( [ -1, 1] \setminus \{ 0 \} ) $ \\

(4.14)~~~ $ ( [ -1, 1] \setminus \{ 0 \} ) \ni r \longmapsto 1 / r \in ( \mathbb{R} \setminus ( -1, 1) ) $ \\

The considerable difference between (4.11), (4.12), and on the other hand, (4.13), (4.14) is obvious. Indeed, in the former two, which describe the self-similar structure of $\mathbb{F}_{\cal U}$, the order reversing bijections represent the set \\

$~~~~~~ mon ( 0 ) \setminus \{ 0 \} $ \\

through the set \\

$~~~~~~ [~ (~ \bigcup_{r \in \mathbb{R}, \,\, \lambda \in \Lambda} \,
                       (r - s_\lambda + monad ( 0 )) ~) $ \\ \\
                        \hspace*{3cm} $   ~\bigcup~ (~ \bigcup_{r \in \mathbb{R}, \,\, \lambda \in \Lambda} \, ( r + s_\lambda + monad ( 0 )) ~) ~] $ \\

which contains {\it uncountably} many translates of the set $mon ( 0 )$. And it is precisely this manifestly {\it rich} self-similarity of the set $mon ( 0 )$ of monads which is the novelty in the {\it non-Archimedean} structure of $\mathbb{F}_{\cal U}$, when compared with the much simpler Archimedean structure of $\mathbb{R}$. This novelty is remarkable since it makes $mon ( 0 )$ have the very same {\it complexity} with the whole of \\

$~~~~~~ \mathbb{F}_{\cal U} \setminus Gal ( 0 ) = [~ (~ \bigcup_{r \in \mathbb{R}, \,\, \lambda \in \Lambda} \,
                       (r - s_\lambda + monad ( 0 )) ~) $ \\ \\
                        \hspace*{5cm} $   ~\bigcup~ (~ \bigcup_{r \in \mathbb{R}, \,\, \lambda \in \Lambda} \, ( r + s_\lambda + monad ( 0 )) ~) ~] $ \\

In this way $mon ( 0 )$, which is but the set of infinitesimals, thus it cannot be represented in terms of the usual field $\mathbb{R}$ of real numbers, turns out to have the very {\it same} complexity as the set $\mathbb{F}_{\cal U} \setminus Gal ( 0 )$ of all infinitely large numbers, which again cannot be represented in terms of the usual field $\mathbb{R}$ of real numbers. \\ \\

{\bf 5. What Possible Relevance in Physics ?} \\

In sharp contrast with the rich self-similar structure of the ultrapower fields $ \mathbb{F}_{\cal U}$, which by necessity are {\it not} Archimedean, and thus exhibit the mentioned wealth of {\it walkable worlds}, all the known theories of Physics are - and have so far been stuck - into only one {\it single} such walkable world, see (3.8), namely, the one given by the usual field $\mathbb{R}$ of real numbers, and the scalars, manifolds, spaces, etc., built upon $\mathbb{R}$. No wonder that so called "infinities in Physics" trouble various theories of Physics, and furthermore, lead to highly questionable ad-hock attempts at solution, such as various so called renormalization methods. \\

On the other hand, whatever quantities prove to be "infinity" in Physics becomes just another usual and regular element in the
ultrapower fields  $ \mathbb{F}_{\cal U}$, thus all algebraic operations can be effectuated upon such quantities, without absolutely any concern or restrictions. \\

And for all that, all it takes is simply to discard the Archimedean Axiom. \\

But beyond the issue of ''infinities'' in Physics lies the far deeper and more important issue of possible {\it different levels of precision} which may be relevant in modelling natural processes. \\

Here, with the acceptance of infinitesimals and infinitely large large scalars, and thus with the necessary setting aside of the Archimedean Axiom, we may obviously face two rather different alternatives, namely

\begin{itemize}

\item the new properties and meanings in reduced power algebras do not correspond to any possible physical meaning, or on the contrary

\item such new properties and meanings which appear in ultrapower fields may possibly correspond to not yet explored physical realities.

\end{itemize}

One such possible such new physical reality is that of

\begin{itemize}

\item Increased, or alternatively, Decreased Precision in Physical Measurements.

\end{itemize}

As a general issue, relating not only to Relativity or the Quanta, but possibly to Classical Physics as well, the
presence of infinitesimal and infinitely large scalars in ultrapower fields may correspond to a new possibility of having no less than two radically different kind of measurements when it comes to their {\it relative precision}. \\

Namely, one has an {\it increased precision} in measurement, when measurement is done in terms of usual finite scalars, and one obtains as result some infinitesimal scalar in such ultrapower fields. \\

Alternatively, the presence of infinitely large scalars in such ultrapower fields may simply indicate that they were obtained in terms of finite scalars, and thus are but the result of a measurement with {\it decreased precision}. \\

In this regard, we can therefore have the following relative situations

\begin{itemize}

\item infinitesimal scalars are the result of increased precision measurements done in terms of finite or infinite scalars,

\item finite scalars are the result of increased precision measurements done in terms of infinite scalars,

\item finite or infinitely large scalars are the result of decreased precision measurements done in terms of infinitesimal scalars,

\item infinitely large scalars are the result of decreased precision measurements done in terms of infinitesimal or finite scalars.

\end{itemize}

And surprisingly, one can also have the following relative situations

\begin{itemize}

\item infinitesimal scalars are the result of increased precision measurements done in terms of some less infinitesimal scalars,

\item infinitesimal scalars are the result of decreased precision measurements done in terms of some more infinitesimal scalars,

\item infinitely large scalars are the result of increased precision measurements done in terms of some more infinitely large scalars,

\item infinitely large scalars are the result of decreased precision measurements
done in terms of some less infinitely large scalars.

\end{itemize}

Indeed, one of the basic features of ultrapower fields is precisely their rich self-similar structure which distinguishes not only between infinitesimal, finite and infinitely large scalars, but also within the infinitely small scalars themselves, as much as within the infinitely large scalars taken all alone. Specifically, infinitesimal
scalars can be infinitely smaller, or on the contrary, infinitely larger than other infinitesimals. And similarly, infinitely large scalars can be infinitely smaller, or on the contrary, infinitely larger than other infinitely large scalars. \\

Here, however, we can note that such a possible interpretation of increased, or decreased precision which is {\it relative}, is in fact not new. Indeed, in terms of usual scalars, be they real or complex, there is a marked dichotomy between finite scalars, and on the other hand, the so called ''infinities'' which may on occasion arise from operations
with finite scalars. And such simple ''formulas'' like $\infty + 1 = \infty$, are in fact expressing that fact. Namely, on one hand, from the point of view of ''infinity'', the finite number 1 has such an increased precision as to be irrelevant with respect to addition, while on the other hand, from the point of view of the finite number 1, the ''infinity'' has such an decreased precision as to be similarly irrelevant when involved in addition. \\

As for the issue of certain universal constants in Physics, given the above possibility in interpretation leading to relative precision measurement, be it increased or decreased, one can reconsider the status of certain universal physical constants, such as for instance, the Planck constant h and the constant c giving the velocity of light in
vacuum. \\
Indeed, when considered from our everyday macroscopic experience, h is supposed to be unusually small, while on the contrary, c is very large. Consequently, one may see h as a sort of ''infinitesimal'', while c then looks like ''infinitely large''. \\

The fact is that, within ultrapower fields, such an alternative view of h and c is possible. Therefore, one may find it appropriate to explore the possible physical meaning, or otherwise, that may possibly be associated with such an interpretation. \\ \\

{\bf Appendix} \\

Given a nonvoid set $\Lambda$, by a {\it filter} ${\cal F}$ on $\Lambda$ one understands a set of subsets $F \subseteq \Lambda$ with the following four properties \\

(A1)~~~ $ {\cal F} \neq \phi $ \\

(A2)~~~ $ \phi \notin {\cal F} $ \\

(A3)~~~ $ F, G \in {\cal F} ~~\Longrightarrow~~ F \cap G \in {\cal F} $ \\

(A4)~~~ $ F \in {\cal F},~ F \subseteq G \subseteq \Lambda ~~\Longrightarrow~~ G \in {\cal F} $ \\

As an example, let $\Lambda = \mathbb{N}$ and let \\

$~~~~~~ {\cal F}_{fre} $ \\

be the set of all subsets $F \subseteq \mathbb{N}$ such that $\mathbb{N} \setminus F$ is a finite subset of $\mathbb{N}$. Then ${\cal F}_{fre}$ is a {\it filter} on $\Lambda = \mathbb{N}$, and it is called the Fr\'{e}chet filter. \\

An {\it ultrafilter} on a nonvoid set $\Lambda$ is any filter ${\cal U}$ which has the additional property \\

(A5)~~~ $ F \subseteq \Lambda ~~~\Longrightarrow~~~ F \in {\cal U} $ ~ or ~ $ \Lambda \setminus F \in {\cal U} $ \\

Clearly, ${\cal F}_{fre}$ is {\it not} an ultrafilter on $\Lambda = \mathbb{N}$, since if we take $F$ as the set of all even numbers, then $\Lambda \setminus F$ is the set of all odd numbers, and none of them belongs to ${\cal F}_{fre}$. \\

The existence of ultrafilters on every nonvoid set $\Lambda$ is guaranteed by the Axiom of Choice in Set Theory. \\

{~~} \\ \\

{\LARGE \bf PART II : \\ \\ Do Computations in Calculus \\ \\ Involve Time-Monads ? \\

Is Space-Time Non-Archimedean ?} \\ \\

{\bf Abstract} \\

Whenever a limit $\lim_{\,n \to \infty} x_n = x $ is computed according to usual Calculus, the respective computation, done either by a human mind, or by
some other entity, may by necessity be {\it unconventional} with respect to the usual assumption on the structure of time as given by the set
$\mathbb{R}$ of real numbers. That is, the computation of such limits may in fact take place in richer and more complicated time structures that contain
non-zero infinitesimals, unlike $\mathbb{R}$ which has 0 as the only infinitesimal. Thus the computation of limits according to Calculus may need unconventional computation which involves time intervals better described by {\it time-monads}. A large class of easy to construct and use such richer and more complicated
time structures which have infinitesimals, and thus monads as well, is that of the so called {\it reduced power fields}, or RPF-s, which are particular
cases of {\it reduced power algebras}, or RPA-s. Such RPA-s, and in particular RPF-s, are by necessity non-Archimedean. It may therefore appear, no matter how surprisingly, that at least since the invention of Calculus by Newton, some of us humans have - without being quite aware of it - been doing lots of unconventional computation whenever we computed limits, and of course, derivatives and integrals, which are essentially based on limits. On the other hand, one may see the presence of unconventional computation in Calculus being rather strongly supported by the well documented difficulty experienced by masses of students when faced with learning that mathematical discipline. \\ \\

{\bf 0. Preliminaries} \\

There is already a well entrenched research field called Unconventional Computation which has progressed, and continues to do so quite far beyond the
classical paradigm of Turing machines, [1]. \\
A main avenue pursued in this regard appears to be that in which time and/or space as usually modelled by $\mathbb{R}$ and $\mathbb{R}^3$, respectively,
are replaced with other, and often more rich and complicated mathematical structures. \\
So far, however, unconventional computation has not been claimed to have been implemented effectively in any way whatsoever, except for what may be called ''on paper'', that is, within the respective theories. \\

Here we shall address the {\it pragmatically foundational} question whether {\it effective} unconventional computations are at all possible within {\it structures of time} which are more more rich and complicated than the usual one modelled by $\mathbb{R}$, that is, structures which contain the set $\mathbb{R}$ of usual time moments as a subset. By {\it effective} computation we mean computations which by some existing and easily accessible ways, be they human minds or some other entities, can be performed within usual finite time intervals. It follows that, in case richer time structures than that modelled by $\mathbb{R}$ may be involved, such effective computations may indeed be unconventional. \\

We are concerned about {\it structures of time} since computation, understood in quite general terms, is supposed - in its non-parallel versions - to mean a linearly ordered sequence of information processing operations. And such a sequence - including the assumed duration of the operations of processing information that happen at each step of the sequence - is of course supposed to take place within an appropriate time structure, be that given as usually by $\mathbb{R}$, or by a richer structure, such as for instance, some RPF. \\

What is suggested here in this regard is that the usual operations of limit, namely \\

(0.1)~~~ $ \lim_{\,n \to \infty} x_n = x $ \\

where $n \in \mathbb{N},~ x_n,~ x \in \mathbb{R}$, of the so called 101 Differential and Integral Calculus - performed mentally, or through physical processes described by Classical Mechanics, or possibly, by other similarly {\it effective} means - do in fact already amount to {\it unconventional computations}, since the time structures in which they happen, and much unlike the usual time, appear by necessity to contain infinitesimals, and thus the celebrated {\it monads} of Leibniz. \\

In this regard it should be recalled that the whole of Differential and Integral Calculus is essentially based on operations of limit of type (0.1), the rest being but reducible to various applications of Elementary, that is, pre-Calculus Mathematics. \\

It follows that, quite unknown to us, and at least ever since the discovery of Calculus, time has in fact had an {\it non-Archimedean} structure for all those who have known, understood and used that mathematical discipline. \\
And from here, to a non-Archimedean structure of space as well, there is but one small further step ... \\ \\

{\bf 1. Computation, Physics and Calculus} \\

There is by now a more than three century long historical record about the following two facts :

\begin{itemize}

\item Classical Mechanics cannot be developed based alone on pre-Calculus type Elementary Mathematics,

\item Classical Mechanics can be developed based on Calculus.

\end{itemize}

Obviously, the above two facts alone are not a proof that, indeed, Calculus is {\it necessary} for the precise mathematical modelling of Classical Mechanics. Instead, they only show the well known fact that it is {\it sufficient} for such a modelling. In other words, it has so far not been proved as impossible to develop a mathematical theory strictly weaker than Calculus, yet still capable to model in a precise mathematical manner Classical Mechanics. \\
This is why we shall focus on Calculus, and not on Classical Mechanics when dealing with unconventional
computation. \\

As for Newton, he obviously created Calculus precisely in view of the fact that, even if tacitly, he considered Calculus as being necessary for a precise mathematical formulation of Classical Mechanics. Fortunately, he succeeded, since Calculus proved to be sufficient for a precise mathematical formulation of Classical Mechanics. \\

Let us start with two simple examples given by well known formulas in Calculus \\

(1.1)~~~ $ 1 + \frac{1}{2 !} + \frac{1}{3 !} + \ldots = e $ \\

(1.2)~~~ $ \lim_{\,n \to \infty} \frac{n}{n + 1} = 1 $ \\

The terms of the infinite series in (1.1), just as those of the infinite sequence in (1.2), can be generated by two respective simple algorithms. However, what appears to be the essential novel phenomenon brought in by Calculus, when compared with the earlier Elementary Mathematics, is that

\begin{itemize}

\item both relations (1.1) and (1.2) mean in a certain sense that infinitely many arithmetic operations are performed within a finite time interval, in order to establish these relations.

\end{itemize}

Furthermore, the result $e$ in (1.1) is a transcendental - and not merely an algebraic - real number, and as such, it contains an infinite amount of information. \\

Therefore, the conclusion is that

\begin{itemize}

\item Calculus offers the possibility to do within a finite usual time interval, with time meant in the usual sense as modelled by the field $\mathbb{R}$ of real numbers, infinitely many arithmetic operations.

\end{itemize}

What is further remarkable in (1.1) and (1.2) is that in such a performance by Calculus there is {\it no} any kind of Zeno-type effect within usual time when it comes to the complexity of computations involved at successive steps. In fact, we may clearly note an opposite effect. Indeed, each new term in the infinite series in (1.1), and similarly, in the infinite sequence in (1.2), leads to a more complex computation than the previous one. Therefore, the above conclusion can be further strengthened as follows

\begin{itemize}

\item Calculus offers the possibility to do within a finite usual time interval infinitely many, and more and more complex arithmetic operations.

\end{itemize}

Related to (1.1) and (1.2) one may argue in the following alternative manner. In series like (1.1) one does actually have a mere {\it definition} of the number $e$, or of a corresponding other one, and as such, no any kind of infinite calculations need be involved. As for sequences like (1.2) one can suggest that the result $1$, or a corresponding other one in the right hand term can be {\it guessed} based on suitable arguments, thus again, no any kind of infinite calculations need be involved. \\

However in case of such an argument one issue remains open, namely, how does one know that the infinite series in (1.1), and the infinite sequence in (1.2) are indeed convergent ? \\

In this way, one is left with having to apply the Cauchy Criterion, in order to establish the respective convergence properties. And then, one cannot avoid getting involved in universal logical quantifiers ranging over infinite domains, an issue addressed in the next section. \\

The effect of the above quite naturally is to ask the questions

\begin{itemize}

\item How can Calculus offer the possibility to do within a finite usual time interval infinitely many, and more and more complex arithmetic operations ?

\item Which is the kind of time structure within which Calculus manages such a performance ?

\end{itemize}

\bigskip

{\bf 2. The Universal Logical Quantifier} \\

As if to further aggravate the situation, there comes the rigorous definition of (0.1) according to Calculus, namely \\

(2.1)~~~ $ \forall~ \epsilon > 0 ~:~ \exists~ m \in \mathbb{N} ~:~ \forall~ n \in \mathbb{N} ~:~ n \geq m
                                                                ~\Longrightarrow~ |\, x - x_n \,| \leq \epsilon $ \\

The essential fact in (2.1) from the point of view of Calculus is that the two universal quantifiers $\forall$ which appear in it range over the infinite
domains $\epsilon \in ( 0, \infty )$ and $n \in \mathbb{N}$, respectively. Yet in the mind of a human being who knows and understands Calculus, the
operations of the respective two universal quantifiers happen within a finite usual time interval. \\

As for the Cauchy Criterion for the convergence of an infinite sequence like in (0.1), we have the rigorous formulation given by \\

(2.2)~~~ $ \begin{array}{l}
             \forall~ \epsilon > 0 ~:~ \exists~ k \in \mathbb{N} ~:~ \forall~ n \in \mathbb{N} ~:~
               \forall~ m \in \mathbb{N} ~:~ \\ \\
                 ~~~~ n, m \geq k ~\Longrightarrow~ |\, x_m - x_n \,| \leq \epsilon
            \end{array} $ \\

thus this time a triple - and not only double like in (2.1) - involvement of the universal quantifier ranging over infinite domains. \\

Furthermore, the problem of dealing in (2.1) and (2.2) with infinity and doing so within a finite usual time interval is obviously not limited only to the universal quantifiers involved, since it also occurs in the existential quantifiers present as they similarly range over infinite domains. \\

But then, such a mental process in humans is obviously a particular case of the ability of human mind to conceive infinity, be it actual or
potential, and do so within a finite usual time interval. \\

The issue of potential versus actual infinity has a long history, and its importance both within, and outside of Mathematics cannot easily be
overstated. \\
In this regard, however, it can be noted that this issue has mostly been addressed from its ontological perspective which imposed its various
conclusions upon possible gnoseologic, epistemic, let alone pragmatic views. \\

One of the major novelties, therefore, brought about by Calculus, and specifically, by its quintessential operation of limit in (0.1), (2.1), is to place
the {\it pragmatic} aspect of infinity up front, and in fact, to assume the human ability to deal with an actually infinite amount of arithmetical
operations, and do so within a finite usual time interval. \\

In this regard, the ancient paradoxes of Zeno - to the extent that they are taken as being genuine paradoxes - appear to be no more than an expression of a mental inability to make the very last step done by Newton,
namely, to jump from the potential infinity in a sequence $x_1, x_2, x_3, \ldots , x_n, \ldots $ to $\lim_{\,n \to \infty} x_n = x$, seen as an actual
infinity. \\

In this way, Calculus, within its specific realms, has given a first and major treatment of the age old issue of potential versus actual infinity, and
has done so pragmatically, and as such, in massively useful ways. \\
However, the most impressive theoretical approach to infinity has, so far, been the Set Theory
of Cantor which, clearly, goes ways beyond Calculus. \\

If we are concerned about the issue of infinity, we then can recall that, together with the issues of self-referentiality and change, it had already been
of major importance in pre-Socratic philosophy, [6]. In this regard, Calculus has brought a pragmatic approach not only to infinity, but also to change.
As for self-referentiality, except for the well known paradoxes it has led ever since ancient times, its most spectacular and consequential modern
approach has been that resulting in G\"{o}del\,'s Incompleteness Theorems. \\

When considering these three issues of self-referentiality, change and infinity, it may not be inappropriate to see them as somehow essentially related
on more subtle levels. A hint in this regard, among many other ancient or more modern ones, may be in Exodus 3:14 in the Old Testament, where the name
of God is given by what is claimed to be the all encompassing and ultimately dynamic self-referentiality of "I AM THAT I AM". \\ \\

{\bf 3. A Class of Enriched Time Structures Large and Easy to \\
        \hspace*{0.45cm} Construct and Use} \\

Abraham Robinson's 1966 construction of Nonstandard Analysis, and specifically, of the nonstandard real line $^*\mathbb{R}$, did in fact introduce, among others, time and space structures significantly more rich and complicated
than the classical ones. And it did so along the well known, albeit mathematically vague, 17th century ideas of Leibniz regarding {\it infinitesimals}, and therefore,
{\it monads}. \\

The fact that Nonstandard Analysis has nevertheless failed to become popular even within the larger mathematical community, let alone within other science disciplines, is due to what can be seen as a sort of cost-return
situation in which the vast majority of so called ''working mathematicians'' have, rightly or wrongly, decided that the returns do in no way justify the costs, when learning the rather complex mathematics involved in Nonstandard Analysis. \\

As it happens, however, the very same enriched time and space structures given by Nostandard Analysis - and called reduced power algebras, or RPA-s - can be obtained in a far simpler manner, namely, using only 101 Algebra, and specifically, the concepts of ring, ideal, quotient, as well as the rather simple and intuitive set theoretic concept of filter. \\
What one loses by that much more simple and easy approach is the Transfer Principle in Nonstandard Analysis. However, as seen in the sequel, such a loss does not inconvenience to any significant extent. \\
Furthermore, as is well known, the Transfer Principle in Nonstandard Analysis suffers from a severe limitation, namely, it is restricted to entities which can be described by what is called First Order Predicate Logic. On the other hand, a large amount of entities in Calculus, not to mention the rest of Mathematics, cannot be formulated within First Order Predicate Logic, thus fall outside of the range of applicability of the Transfer Principle. \\

We recall here briefly the construction of the so called reduced power algebras, or RPA-s, which have been well known in Model Theory, where they constitute one of the most important basic concepts. Notations and details used in the sequel can be found in [2-5]. \\

This construction happens in the following {\it three} steps

\begin{itemize}

\item first, one chooses an arbitrary {\it infinite} index set $\Lambda$ and constructs the power algebra $\mathbb{R}^\Lambda$ which is but the set of all functions $f : \Lambda \longrightarrow \mathbb{R}$ considered with the usual addition and multiplication of functions,

\item second, one chooses a {\it proper} ideal ${\cal I}$ in the power algebra $\mathbb{R}^\Lambda$,

\item third, one constructs the quotient algebra

\end{itemize}

(3.1)~~~ $ A = \mathbb{R}^\Lambda / {\cal I} $ \\

which is called a {\it reduced power algebra}, or RPA. \\

An important simplification of this construction can be obtained by replacing proper ideals ${\cal I}$ in the algebra $\mathbb{R}^\Lambda$, with the
simpler mathematical structures of {\it filters} ${\cal F}$ on $\Lambda$. Here we recall that a {\it filter} ${\cal F}$ on $\Lambda$ is a set of subsets
$I \subseteq \Lambda$ with the following four properties \\

(3.2)~~~ $ {\cal F} \neq \phi $ \\

(3.3)~~~ $ \phi \notin {\cal F} $ \\

(3.4)~~~ $ I, J \in {\cal F} \Longrightarrow I \cap J \in {\cal F} $ \\

(3.5)~~~ $ I \in {\cal F},~~ I \subseteq J \subseteq \Lambda \Longrightarrow J \in {\cal F} $ \\

Thus such filters can be seen as collections of {\it large} subsets of $\Lambda$. Indeed, (3.2) means that there exist such large subsets, and certainly,
none of them is void, as required by (3.3). Condition (3.4) means that the intersection of two large subsets is still a large subset, while
(3.5) simply means that a subset which contains a large subset is itself large. In particular, $\Lambda$ itself is large, thus $\Lambda \in {\cal F}$. \\

And important example of filter on $\Lambda$ is the {\it Frech\'{e}t} filter, given by \\

(3.6)~~~ $ {\cal F}_{re} ( \Lambda ) = \{~ I \subseteq \Lambda ~~|~~ \Lambda \setminus I ~~\mbox{is finite} ~\} $ \\

The mentioned simplification comes about through the following one-to-one simple correspondence between proper ideals ${\cal I}$ in the algebra
$\mathbb{R}^\Lambda$ and filters ${\cal F}$ on $\Lambda$, namely \\

(3.7)~~~ $ {\cal I} ~~\longmapsto~~ {\cal F}_{\cal I} = \{ Z ( x ) ~|~ x \in {\cal I} \} ~~\longmapsto~~ {\cal I}_{{\cal F}_{\cal I}} = {\cal I} $ \\

(3.8)~~~ $ {\cal F} ~~\longmapsto~~ {\cal I}_{\cal F} = \{ x \in \mathbb{R}^\Lambda ~|~ Z ( x ) \in {\cal F} \}
                         ~~\longmapsto~~ {\cal F}_{{\cal I}_{\cal F}} = {\cal F}$ \\

where for $x \in \mathbb{R}^\Lambda$ we denote $Z ( x ) = \{ \lambda \in \Lambda ~|~ x ( \lambda ) = 0 \}$. \\

An important consequence of (3.1), (3.7), (3.8) is that the mapping \\

(3.9)~~~ $ \mathbb{R} \ni r \longmapsto u_r + {\cal I} \in A = \mathbb{R}^\Lambda / {\cal I} $ \\

is an {\it injective} algebra homomorphism for every proper ideal ${\cal I}$ in $\mathbb{R}^\Lambda$, where $u_r \in \mathbb{R}^\Lambda$ is defined by
$u_r ( \lambda ) = r$, for $\lambda \in \Lambda$. Indeed, in view of (3.7), (3.8), we have for $r \in \mathbb{R}$ \\

(3.10)~~~ $ u_r \in {\cal I} ~~\Longrightarrow~~ r = 0 $ \\

since ${\cal I} = {\cal I}_{{\cal F}_{\cal I}}$, while $u_r \in {\cal I}_{{\cal F}_{\cal I}}$ gives $Z ( u_r ) \in {\cal F}_{\cal I}$, thus $Z ( u_r )
\neq \phi$, which means $r = 0$. \\

Of interest here are a particular case of filters on $\Lambda$, called {\it ultrafilters} ${\cal U}$, and which are characterized by the property \\

(3.11)~~~ $ \forall~~ I \subseteq \Lambda ~:~ I \notin {\cal U} \Longrightarrow \Lambda \setminus I \in {\cal U} $ \\

One of their properties relevant in the sequel is that, through (3.7), (3.8), ultrafilters are in one-to-one correspondence with {\it maximal} ideals in
$\mathbb{R}^\Lambda$, namely \\

(3.12)~~~ $ {\cal F} ~\mbox{ultrafilter} ~~\Longrightarrow~~ {\cal I}_{\cal F} ~\mbox{maximal ideal} $ \\

(3.13)~~~ $ {\cal I} ~\mbox{maximal ideal} ~~\Longrightarrow~~ {\cal F}_{\cal I} ~\mbox{ultrafilter} $ \\

For our purposes, it is useful to distinguish between {\it fixed}, and on the other hand, {\it free} ultrafilters on $\Lambda$. The fixed ones are of the
form ${\cal U}_\lambda = \{ I \subseteq \Lambda ~|~ \lambda \in I \}$, for any given $\lambda \in \Lambda$, while the free ones are all the other
ultrafilters on $\Lambda$. It is easy to see that an ultrafilter ${\cal U}$ on $\Lambda$ if free, if and only if \\

(3.14)~~~ $ {\cal F}_{re} ( \Lambda ) \subseteq {\cal U} $ \\

and the existence of free ultrafilters results from the Axiom of Choice. \\

Now we recall from Algebra that \\

(3.15)~~~ $ {\cal I} ~\mbox{maximal ideal in}~ \mathbb{R}^\Lambda ~~\Longleftrightarrow~~ \mathbb{R}^\Lambda / {\cal I} ~\mbox{field} $ \\

And then (3.12) - (3.15) will result in \\

{\bf Theorem 3.1.} \\

Let ${\cal U}$ be a filter on $\Lambda$ for which (3.14) holds. Then \\

(3.16)~~~ $ \mathbb{F}_{\cal U} = \mathbb{R}^\Lambda / {\cal I}_{\cal U} $~ is a field
                      ~~$ \Longleftrightarrow~~ {\cal U} $~ is an ultrafilter on ~$ \Lambda $ \\

in which case \\

(3.17)~~~ $ \mathbb{R} \subsetneqq \mathbb{F}_{\cal U} $ \\

with the corresponding injective and non-surjective algebra homomorphism, see (3.9) \\

(3.18)~~~ $ \mathbb{R} \ni t \longmapsto u_t + {\cal I}_{\cal U} \in \mathbb{F}_{\cal U} =
                                                         \mathbb{R}^\Lambda / {\cal I}_{\cal U} $ \\

and \\

(3.19)~~~ $ \mathbb{F}_{\cal U} = \mathbb{R}^\Lambda / {\cal I}_{\cal U} $ \\

is called an {\it ultrapower field}, or in short, UPF.

\hfill $\Box$ \\

The important novelty with UPF-s such as $\mathbb{F}_{\cal U}$ in (3.19) is that they are constituted from the following {\it three} different kind of elements $t \in \mathbb{F}_{\cal U}$, called respectively {\it infinitesimal, finite} and {\it infinite}, namely \\

(3.20)~~~ $ \forall~~ r \in \mathbb{R},~ r > 0 ~:~ t \in ( - r, r ) $ \\

(3.21)~~~ $ \exists~~ r \in \mathbb{R},~ r > 0 ~:~ t \in ( - r, r ) $ \\

(3.22)~~~ $ \forall~~ r \in \mathbb{R},~ r > 0 ~:~ t \notin ( - r,  r ) $ \\

where for $a, b \in \mathbb{F}_{\cal U}$, we denote as usual $( a, b ) = \{ s \in \mathbb{F}_{\cal U} ~|~ a < s < b \}$. \\

Now,
following Leibniz, one denotes \\

(3.23)~~~ $ monad ( 0 ) = \{~ t \in \mathbb{F}_{\cal U} ~|~~ t ~~\mbox{is~infinitesimal}~ \} $ \\

and calls it the {\it monad} of $0 \in \mathbb{F}_{\cal U}$, while following Keisler, [7], one denotes \\

(3.24)~~~ $ Gal ( 0 ) = \{~ t \in \mathbb{F}_{\cal U} ~|~~ t ~~\mbox{is~finite}~ \} $ \\

and calls it the {\it Galaxy} of $0 \in \mathbb{F}_{\cal U}$. \\

It is easy to see that \\

(3.25)~~~ $ Gal ( 0 ) = \bigcup_{r \in \mathbb{R}} monad ( r ) $ \\

where for $t \in \mathbb{F}_{\cal U}$, we denote \\

(3.26)~~~ $ monad ( t ) = t + monad ( 0 ) $ \\

Finally \\

(3.27)~~~ $ \mathbb{F}_{\cal U} \setminus Gal ( 0 ) $ \\

is the set of infinitely large elements in the ultrapower field $\mathbb{F}_{\cal U}$. \\

In this way, all the elements of $\mathbb{F}_{\cal U}$, be they infinitesimal, finite, or infinitely large, have been expressed respectively in (3.23) by the monad of $0 \in \mathbb{F}_{\cal U}$, in (3.25) by the Galaxy of $0 \in \mathbb{F}_{\cal U}$, and finally, in (3.27). \\

And as one notes, all these sets can in fact be expressed in terms of the monad of $0 \in \mathbb{F}_{\cal U}$ alone. \\ \\

{\bf 4. Rich and Complicated Self-Similar Structure of Reduced
        \hspace*{0.55cm} Power Fields} \\

For the sake of easier comparison, let us recall the {\it self-similarity} properties of the usual field $\mathbb{R}$ of real numbers. In this regard, we have the {\it self-similarity} property given by the following {\it bijective, order reversing} mapping \\

(4.1)~~~ $ \mathbb{R} \setminus ( -1, 1 ) \ni r \longmapsto 1 / r \in [ -1, 1 ] \setminus \{ 0 \} $ \\

thus the unbounded set \\

$~~~~~~ \mathbb{R} \setminus ( -1, 1 ) = ( -\infty, -1 ] \cup [ 1, \infty ) $ \\

has through the mapping (4.1) the inverse linear order structure of the bounded set \\

$~~~~~~ [ -1, 1 ] \setminus \{ 0 \} = [ -1, 0 ) \cup ( 0, 1 ] $ \\

Now by translation and scaling, we obtain the family of self-similarities of the usual field $\mathbb{R}$ of
 real numbers, given by the bijective, order reversing mappings \\

(4.2)~~~ $ \mathbb{R} \setminus ( -a, a ) \ni r \longmapsto
                         ( 1 / r ) + r_0 \in [ r_0 - \frac{1}{a}, r_0 + \frac{1}{a} ] \setminus \{ r_0 \} $ \\

where $r_0, a \in \mathbb{R}, a > 0$. \\

Here we can note that none of the self-similarities (4.2) refers to the structure at any given point $r_0 \in \mathbb{R}$, but only to the structure of the sets \\

(4.3)~~~ $ [ r_0 - a, r_0 + a ] \setminus \{ r_0 \} = [ r_0 - a, r_0 ) \cup ( r_0, r_0 + a ],~~~~ a > 0 $ \\

around points $r_0 \in \mathbb{R}$, sets which are whole neighbourhoods of $r_0$ from which, however, the point $r_0$ itself has been taken out. This is obviously inevitable, since each point $r_0 \in \mathbb{R}$ is at a finite strictly positive - thus not infinitesimal - distance from any other point in $\mathbb{R}$. \\

In addition, we also have the self-similarities \\

(4.4)~~~ $ \mathbb{R} \stackrel{f}\longrightarrow ( a, b ) $ \\

where $-\infty \leq a < b \leq \infty$, while $f$ can be any bijective order preserving continuous mapping. \\

On the other hand, with the ultrapower fields  $ \mathbb{F}_{\cal U}$, their self-similar structures are far more {\it
rich}, due to the presence of their infinitesimals, and thus as well, of their infinitely large elements. Indeed, this
time, the self-similarities can also refer to the whole {\it monad} of each point, except for the point itself. \\

Let us start with a self-similarity of any ultrapower field $ \mathbb{F}_{\cal U}$ which does {\it not} exist in the
case of the usual real line $\mathbb{R}$. Namely, it is easy to see that we have the order reversing bijective mapping \\

(4.5)~~~ $ ( \mathbb{F}_{\cal U} \setminus Gal(0) ) \ni t \longmapsto 1 / t \in ( monad(0) \setminus \{ 0 \} ) $ \\

which means that the set of all infinitely large elements in $\mathbb{F}_{\cal U}$ has the inverse order structure of the set of infinitesimal elements from which one excludes 0. \\

This shows the {\it important} fact that the infinitesimally local structure, and on the other hand, the global structure of $\mathbb{F}_{\cal U}$ do in fact {\it mirror} one another, a property which has no correspondence in the case of the usual field $\mathbb{R}$ of real numbers. \\

Also, through translation and scaling, we have, for each $t_0, u \in \mathbb{F}_{\cal U}, u > 0$, the order reversing bijective mapping \\

(4.6)~~~ $  \mathbb{F}_{\cal U} \setminus ( -u, u ) \ni t \longmapsto
                      ( 1 / t ) + t_0 \in [ t_0 - \frac{1}{u}, t_0 + \frac{1}{u} ] \setminus \{ t_0 \} $ \\

where $\mathbb{F}_{\cal U} \setminus ( -u, u )$ will always contain infinitely large elements. \\

These again are self-similarities not present in the case of the usual real line $\mathbb{R}$. \\

Furthermore, in (4.6) we have a far more rich possibility for translations and scalings than in the usual case of the real line $\mathbb{R}$. Indeed, in addition to translations and scalings with non-zero finite elements $r_0, a \in \mathbb{R}, a > 0$, as in (4.3), we can now also translate and scale with all $t_0, u \in \mathbb{F}_{\cal U}, u > 0$, thus with all infinitely large elements, as well as with all infinitesimal elements, except for scaling with $0 \in \mathbb{F}_{\cal U}$. \\

Let us consider the above in some detail by listing the different possibilities for the sets \\

(4.7)~~~ $ [ t_0 - \frac{1}{u}, t_0 + \frac{1}{u} ] \setminus \{ t_0 \} $ \\

in (4.6). \\

First of all, these sets are no longer mere subsets in $\mathbb{R}$, but instead, they are subsets in $\mathbb{F}_{\cal U}$, and will always contain infinitesimals, since they contain nonvoid intervals. Furthermore, as seen below, they may also contain infinitely large elements. \\

Also, $t_0, u \in \mathbb{F}_{\cal U}, u > 0$ in (4.7) can independently be finite, infinitesimal, or infinitely
large, thus resulting in 9 possible combinations and 6 distinct outcomes regarding the set (4.7), which we list below. This is in sharp contradistinction with the case in (4.3) which applies to the real line $\mathbb{R}$. Indeed : \\

1) Let us start the listing of these 9 different cases and 6 distinct outcomes with both $t_0$ and $u$ being finite. Then obviously (4.7) is a subset of $Gal(0)$, and it has the finite, non-infinitesimal length $2u$. \\

2) When $t_0$ is finite and $u$ is infinitesimal, then the set (4.7) is infinitely large, and is no longer contained in $Gal ( 0 )$, however, it contains $Gal ( 0 ) \setminus \{ t_0 \}$. \\

3) If $t_0$ is finite, but $u$ is infinitely large, then (4.7) is again a subset of $Gal ( 0 )$, and in fact, it has the infinitesimal length $2 u$, which means that it is a subset of $monad ( t_0 )$. \\

4) Let us now assume that $t_0$ is infinitesimal and $u$ finite. Then regarding the set (4.1), we are back to case 1) above. \\

5) If both $t_0$ and $u$ are infinitesimal then the set (4.7) is as in 2) above. \\

6) When $t_0$ is infinitesimal and $u$ is infinitely large, the set (4.7) is as in 3) above. \\

7) Let us now take $t_0$ infinitely large and $u$ finite. Then the set (4.7) is disjoint from $Gal ( 0 )$, and it has the finite, non-infinitesimal length $2 u$. \\

8) When $t_0$ infinitely large and $u$ infinitesimal, then the set (4.7) is again not contained in $Gal ( 0 )$, and it has the infinitely large length $2 u$. Furthermore, depending on the relationship between $|t_0|$ and $1 / u$, it may, or it may not intersect $Gal ( 0 )$. \\

9) Finally, when both $t_0$ and $u$ are infinitely large, then the set (4.7) is disjoint from $Gal ( 0 )$, and it has the infinitesimal length $2 u$. \\

We conclude that he {\it local} structure of $\mathbb{F}_{\cal U}$ is given by \\

(4.8)~~~ $ Gal ( 0 ) = \bigcup_{r \in \mathbb{R}} \, ( r + monad ( 0 ) ) $ \\

while the {\it global} structure of $\mathbb{F}_{\cal U}$ is given by \\

(4.9)~~~ $ \mathbb{F}_{\cal U} = (~ \bigcup_{\lambda \in \Lambda} \, Gal ( -s_\lambda ) ~) ~\bigcup~ Gal ( 0 )
                                              ~\bigcup~ (~ \bigcup_{\lambda \in \Lambda} \, Gal ( s_\lambda ) ~) $ \\

where $\Lambda$ is an uncountable set of indices, while $s_\lambda \in \mathbb{F}_{\cal U}$ are positive infinite, and such that $s_\mu - s_\lambda$ is infinite, for $\lambda, \mu \in \Lambda, \lambda \neq \mu$. \\

Here we can point to a self-similar aspect of the interrelation between the local and global structure of $\mathbb{F}_{\cal U}$ which may remind us of a typical feature of {\it fractals}. Indeed, similar with (4.8), the relation (4.9) can also be expressed in terms monads, namely \\

(4.10)~~~ $ \mathbb{F}_{\cal U} = (~ \bigcup_{r \in \mathbb{R}, \,\, \lambda \in \Lambda} \,
                           (r - s_\lambda + monad ( 0 )) ~) ~\bigcup~ $ \\ \\
          \hspace*{3cm} $ (~ \bigcup_{r \in \mathbb{R}} \, ( r + monad ( 0 ) ) ~) $ \\ \\
          \hspace*{3cm} $~\bigcup~ (~ \bigcup_{r \in \mathbb{R}, \,\, \lambda \in \Lambda} \, ( r + s_\lambda +
                            monad ( 0 )) ~) $ \\

In this way, in view of (4.5), we obtain the self-similar order reversing bijection, which is now expressed solely in terms of $mon ( 0 )$, namely \\

(4.11)~~~ $ [~ (~ \bigcup_{r \in \mathbb{R}, \,\, \lambda \in \Lambda} \,
                       (r - s_\lambda + monad ( 0 )) ~) $ \\ \\
                        \hspace*{2cm} $   ~\bigcup~ (~ \bigcup_{r \in \mathbb{R}, \,\, \lambda \in \Lambda} \, ( r + s_\lambda + monad ( 0 )) ~) ~] \ni t \longmapsto $ \\ \\
                        \hspace*{7cm} $ \longmapsto 1 / t \in [~ monad  ( 0 ) \setminus \{ 0 \} ~] $ \\

and conversely \\

(4.12)~~~ $ [~ monad  ( 0 ) \setminus \{ 0 \} ~] \ni t \longmapsto $ \\ \\
               \hspace*{5cm} $ \longmapsto 1/ t \in [~ (~ \bigcup_{r \in \mathbb{R}, \,\, \lambda \in \Lambda} \,
                       (r - s_\lambda + monad ( 0 )) ~) $ \\ \\
                        \hspace*{7cm} $   ~\bigcup~ (~ \bigcup_{r \in \mathbb{R}, \,\, \lambda \in \Lambda} \, ( r + s_\lambda + monad ( 0 )) ~) ~] $ \\

As we can note, the above bijections in (4.11), (4.12) are given by the very simple algebraic, explicit, and order reversing mapping $s \longmapsto 1 / s$, which involves what is essentially a {\it field} operation, namely, division.
And these two bijections take the place of the much simpler order reversing bijections in the case of the usual real line $\mathbb{R}$, namely \\

(4.13)~~~ $ ( \mathbb{R} \setminus ( -1, 1) ) \ni r \longmapsto 1 / r \in  ( [ -1, 1] \setminus \{ 0 \} ) $ \\

(4.14)~~~ $ ( [ -1, 1] \setminus \{ 0 \} ) \ni r \longmapsto 1 / r \in ( \mathbb{R} \setminus ( -1, 1) ) $ \\

The considerable difference between (4.11), (4.12), and on the other hand, (4.13), (4.14) is obvious. Indeed, in the former two, which describe the self-similar structure of $\mathbb{F}_{\cal U}$, the order reversing bijections represent the set \\

$~~~~~~ mon ( 0 ) \setminus \{ 0 \} $ \\

through the set \\

$~~~~~~ [~ (~ \bigcup_{r \in \mathbb{R}, \,\, \lambda \in \Lambda} \,
                       (r - s_\lambda + monad ( 0 )) ~) $ \\ \\
                        \hspace*{3cm} $   ~\bigcup~ (~ \bigcup_{r \in \mathbb{R}, \,\, \lambda \in \Lambda} \, ( r + s_\lambda + monad ( 0 )) ~) ~] $ \\

which contains {\it uncountably} many translates of the set $mon ( 0 )$. And it is precisely this manifestly {\it rich} self-similarity of the set $mon ( 0 )$ of monads which is the novelty in the {\it non-Archimedean} structure of $\mathbb{F}_{\cal U}$, when compared with the much simpler Archimedean structure of $\mathbb{R}$. This novelty is remarkable since it makes $mon ( 0 )$ have the very same {\it complexity} with the whole of \\

$~~~~~~ \mathbb{F}_{\cal U} \setminus Gal ( 0 ) = [~ (~ \bigcup_{r \in \mathbb{R}, \,\, \lambda \in \Lambda} \,
                       (r - s_\lambda + monad ( 0 )) ~) $ \\ \\
                        \hspace*{5cm} $   ~\bigcup~ (~ \bigcup_{r \in \mathbb{R}, \,\, \lambda \in \Lambda} \, ( r + s_\lambda + monad ( 0 )) ~) ~] $ \\

In this way $mon ( 0 )$, which is but the set of infinitesimals, thus it cannot be represented in terms of the usual field $\mathbb{R}$ of real numbers, turns out to have the very {\it same} complexity as the set $\mathbb{F}_{\cal U} \setminus Gal ( 0 )$ of all infinitely large numbers, which again cannot be represented in terms of the usual field $\mathbb{R}$ of real numbers. \\ \\

{\bf 5. One Way to Interpret the Operation of Limit in Calculus} \\

Let us recall the two questions at the end of section 1, namely

\begin{itemize}

\item How can Calculus offer the possibility to do within a finite usual time interval infinitely many, and more and more complex arithmetic operations ?

\item Which is the kind of time structure within which Calculus manages such a performance ?

\end{itemize}

We can start with what may seem to be the first, basic and simplest act when dealing with infinity in Mathematics, be it potential or actual, namely, when we define - and assume to comprehend as well - the set $\mathbb{N}$ of natural numbers. \\
Obviously $\mathbb{N}$ is involved in limits (0.1), therefore infinity in these limits is involved - as a definition or act of comprehension - in no simpler manner than in $\mathbb{N}$ itself. \\

When it comes to limits such as in (2.1) and (2.2), computational aspects involving infinity can also be considered. \\

It follows that we may distinguish between three aspects related to infinity when dealing with limits in (0.1), namely

\begin{itemize}

\item definition

\item comprehension

\item computation

\end{itemize}

Now, the definition aspect, when considered in itself alone, may be seen as being able to avoid issues related to infinity,
since it can be seen as subjected to the only requirement of being logically sound. \\

Therefore, the above two questions may be seen as rather relating to comprehension and computation alone. \\

Regarding computation, it is clear that a Turing machine is not sufficient to deal with (0.1) in general, and not even
with particular cases such as in (2.1) or (2.2). \\

As for comprehension by human minds which deal with Calculus, the fact that such a comprehension can take place - if it happens - within a finite usual time interval is quite clear. \\

However, here we are not so much dealing with human comprehension in general, as rather with the ways Calculus deals with infinity. And then {\it computation} in tandem with comprehension becomes the relevant aspect in the above two questions. \\

When, on the other hand, dealing with human comprehension in general, a whole range of fundamental issues may arise, most of them not yet considered to any satisfactory extent, among them those mentioned in [8]. \\

In the above interactions between computation and comprehension, the mentioned failure of Turing machines, as well as the lack of any Zeno-type effect involved when dealing with such simple examples of limits as in (2.1) or (2.2), for instance, can be seen in view of (4.9) - (4.12) as indicating the possible {\it presence} of time-monads.


\begin{thebibliography}{99}

\bibitem{} Keisler J : Elementary Calculus: An Approach Using Infinitesimals. On-line Edition. \\
http://www.math.wisc.edu/~keisler/calc.html

\bibitem{} Rosinger E E : What scalars should we use ? arXiv:math/0505336

\bibitem{} Rosinger E E : Solving Problems in Scalar Algebras of Reduced Powers
arXiv:math/0508471

\bibitem{} Rosinger E E : Covariance and Frames of Reference. arXiv:quant-ph/0511112

\bibitem{} Rosinger E E : From Reference Frame Relativity to Relativity of Mathematical Models :
Relativity Formulas in a Variety of Non-Archimedean Setups. arXiv:physics/0701117

\bibitem{} Rosinger E E : Cosmic Contact : To Be, or Not To Be \\ Archimedean ? arXiv:physics/0702206

\bibitem{} Rosinger E E : String Theory: a mere prelude to \\ non-Archimedean Space-Time Structures? \\
arXiv:physics/0703154

\bibitem{} Rosinger E E : Mathematics and "The Trouble with Physics", How Deep We Have to Go ?
arXiv:0707.1163

\bibitem{} Rosinger E E : How Far Should the Principle of Relativity Go ? arXiv:0710.0226

\bibitem{} Rosinger E E : Two Essays on the Archimedean versus Non-Archimedean Debate. arXiv:0809.4509

\bibitem{} Rosinger E E : Heisenberg Uncertainty in Reduced Power Algebras. arXiv:0901.4825

\bibitem{} Rosinger E E : No-Cloning in Reduced Power Algebras. \\ arXiv:0902.0264

\bibitem{} Rosinger E E : Special Relativity in Reduced Power Algebras. arXiv:0903.0296

\bibitem{} Rosinger E E : Where Infinitesimals Come From ... \\ arXiv:0909.4396

\end{thebibliography}

\begin{thebibliography}{99}

\bibitem{} UC 2008, 7-th International Conference on Unconventional Computation, August 25-28, Vienna, Austria, Parkhotel Sch\"{o}nbrunn : Physics and
Computation (PC), Optical Supercomputing (OSC), Computing with Biomolecules (CBM)

\bibitem{} Rosinger E E : How Far Should the Principle of Relativity Go ? arXiv:0710.0226

\bibitem{} Rosinger E E : Mathematics and The Trouble with Physics, How Deep We Have to Go ? arXiv:0707.1163

\bibitem{} Rosinger E E : String Theory: a mere prelude to \\ non-Archimedean Space-Time Structures? \\ arXiv:physics/0703154

\bibitem{} Rosinger E E : Real Scalars with a Rich Structure of "Walkable Worlds" : Universes next to Universes, and/or within Universes ... and so on, ad infinitum ... arXiv:0911.4824

\bibitem{} Mortensen C : Change. http://plato.stanford.edu

\bibitem{} Keisler J : Elementary Calculus: An Approach Using Infinitesimals. On-line Edition. \\
http://www.math.wisc.edu/~keisler/calc.html

\bibitem{} Rosinger E E : Where and how does it happen ? \\ arXiv:physics/0505041

\end{thebibliography}
\end{document}